\documentstyle[12pt]{amsart}

\newtheorem{Thm}{Theorem}[section]

\newtheorem{Cor}[Thm]{Corollary}
\newtheorem{Lem}[Thm]{Lemma}
\newtheorem{Prop}[Thm]{Proposition}
\newtheorem{Def}[Thm]{Definition}
\newtheorem{Rmk}[Thm]{Remark}
\newtheorem{Obs}[Thm]{Observation}

\newcommand{\im}{\mbox{$\Rightarrow$}}

\newcommand{\co}{{\rm co}}

\def \N{\rm {\bf N}}
\def \R{\rm {\bf R}}

\begin{document}

\title{Twisted sums, Fenchel-Orlicz spaces and property (M)}
\author{G. Androulakis}
\author{C. D. Cazacu}
\author{N. J. Kalton}
\address{Department of Mathematics \\
University of Missouri-Columbia \\
Columbia, MO 65211
}

\email[G. Androulakis]{giorgis@@math.missouri.edu}
\email[C. D. Cazacu]{cazacu@@picard.math.missouri.edu}
\email[N. J. Kalton]{nigel@@math.missouri.edu}


\maketitle

\noindent
{\bf Abstract:} We study certain twisted sums of Orlicz spaces
with non-trivial type which can be viewed as Fenchel-Orlicz spaces
on $\R^2$. We then show that a large class of Fenchel-Orlicz spaces
on $\R^n$ can be renormed to have property (M). In particular this
gives a new construction of the twisted Hilbert space $Z_2$ and
shows it has property (M), after an appropriate renorming.

\bigskip
\section{Introduction} \label{S:intro}

A twisted sum $Z$ of two Banach spaces $X$ and $Y$ is defined
(see \cite{KP}) through a short exact sequence:
$0\longrightarrow X \longrightarrow Z \longrightarrow Y \longrightarrow 0$.
These short exact sequences in the category of (quasi-)Banach spaces are
considered naturally in the investigation of three space properties
(a property $P$ in the category of quasi-Banach spaces is called
a three space property if for every short exact sequence as above, $Z$
has property $P$ whenever $X$ and $Y$ have it). The roots of this theory
go to Enflo, Lindenstrauss and Pisier's solution \cite {ELP} to Palais'
problem: the property of being isomorphic to a Hilbert space is not a
three space property. The first systematic study of twisted sums of
quasi-Banach spaces appears in \cite{KP}. In that paper
twisted sums of quasi-Banach spaces $X$ and $Y$ are associated to quasi-
linear maps from $Y$ to $X$ and the Banach spaces $Z_p, 1<p<\infty$, are
studied as examples of twisted sums of $\ell_p$'s. In particular,
$Z_2$ is a reflexive Banach space with a basis which has a closed
subspace $X$ isometric to $\ell_2$ with $Z_2/X$ also isometric to $\ell_2$.
$Z_2$ is isomorphic to its dual, yet  $Z_2$ is not isomorphic to $\ell_2$.
Furthermore, $Z_2$ has no
complemented subspace with an unconditional basis, in particular
it has no complemented subspace isomorphic to $\ell_2$. $Z_2$ has an
unconditional finite dimensional Schauder decomposition into two dimensional
spaces (2-UFDD), yet Johnson, Lindenstrauss and Schechtman \cite{JLS} showed
that it fails to have local unconditional structure (l.u.st.).
Twisted sums appear also in a natural way in complex interpolation \cite{K3}.
There are several open problems on twisted sums and in particular on $Z_2$,
(see \cite{K2}), which make the study of these spaces very
interesting.

We will use the class of Fenchel-Orlicz spaces. These spaces were
introduced by Turett \cite{T} and they form a natural generalization
of Orlicz spaces. A main difference between Orlicz spaces and
Fenchel-Orlicz spaces is the replacement of the Orlicz function
defined on $\R_+$
by a Young's function defined on a given normed linear space. The
elements of a Fenchel-Orlicz sequence space will then be sequences in
the given normed linear space.

Property (M) was introduced in \cite{K4} as a tool in the study
of M-ideals of compact operators. In that paper it is proved that for
a separable Banach space $X$, the compact operators form an M-ideal in
the space of bounded operators if and only if $X$ has property (M) and
there is a sequence of compact operators $K_n$ such that $K_n
\rightarrow I$ strongly, $K_n^* \rightarrow I$ strongly and
$\lim_{n\rightarrow \infty} \| I - 2K_n \| =1$. For a detailed study
of M-ideals we refer to \cite{HWW}.

We now give a brief overview of the paper. In Section \ref{S:quasi} we
introduce quasi-convex functions. A function on $\R^n$ is quasi-convex
if and only if it is equivalent to a convex function (Proposition
\ref{P:quasi}). We construct a large class of examples of quasi-convex
maps on $\R^2$(Theorem \ref{T:ex}). In Section \ref{S:TSFO} we show how
quasi-convex maps can replace Young's functions in
generating Fenchel-Orlicz spaces on $\R^n$. The main result of the
section is that a twisted sum of an Orlicz space with type $p>1$  with
itself can be represented as a Fenchel-Orlicz space over $\R^2$
(Theorem \ref{T:twistOFO}). This
includes the case of the spaces $Z_p , 1 < p < \infty$. In Section
\ref{S:FOM} we use a method of \cite{K4} to prove that if $\phi$ is a
Young's function on $\R^n$ which is $0$ only at $0$, then the
Fenchel-Orlicz space $h_\phi$ has property (M) (Theorem
\ref{T:FOproM}). Combining results of the last two sections we see
that the spaces $Z_p , 1 < p < \infty$, have property (M).

\section{Quasi-convex maps} \label{S:quasi}
Let $\R_+$ ( respectively $\overline{\R}_+$) denote the set of non-negative
(respectively extended) real numbers.
\begin{Def}
A function $\phi :\R^n \rightarrow \R_+$ is quasi-convex if there
exists $L>0$ such that  for every $ x_1 , x_2 \in \R^n$ and for every
$\lambda_1 , \lambda_2
\in [0,1]$ with $\lambda_1 + \lambda_2 =1$ we have
\[\phi (\lambda_1 x_1 + \lambda_2 x_2 ) \leq L (\lambda_1 \phi(x_1) +
\lambda_2 \phi(x_2) ) \]
\end{Def}
Note that quasi-convex maps can be defined on a
vector space and that quasi-norms are quasi-convex. In order to give a
characterization of quasi-convex maps on $\R^n$ we introduce an
equivalence relation, standard in the study of Orlicz spaces
(cf. \cite{LT1}).
\begin{Def}
Two functions $\phi$ and $\psi : \R^n \rightarrow \R_+$ are equivalent
($\phi \sim \psi$) if there exists $M>0$ such that $\frac{1}{M}\phi(x) \leq
\psi(x) \leq M \phi(x)$ for all $x\in \R^n$.
\end{Def}
We shall say that two functions are equivalent on a set $B$ if the
above inequalities hold for all $x \in B$. We recall that the convex
envelope of a function $\phi:\R^n \rightarrow \R_+$ is defined by:
\begin{eqnarray}
\co\phi(t)\stackrel{\rm def}{=}\inf \{ \sum_i \alpha_i \phi(t_i)
: t=\sum_i \alpha_i t_i , \mbox{ where } t_i \in \R^n, \sum_i \alpha_i
=1 , \alpha_i \geq 0\}. \nonumber
\end{eqnarray}
It is easy to see that
\begin{itemize}
\item co$\phi(t) \leq \phi(t)$ for all $t \in \R^n$ and
\item if $\psi : \R^n \longrightarrow \R_+$ is a convex
function with $\psi(t) \leq \phi(t)$ for all $t \in \R^n$, then
$\psi(t) \leq$co$\phi(t)$ for all $t \in \R^n$.
\end{itemize}
\begin{Prop} \label{P:quasi}
Let $\phi: \R^n \longrightarrow \R_+$. The following are equivalent:
\begin{itemize}
\begin{enumerate}
\item{ $\phi$ is quasi-convex.}
\item{ $\phi \sim \co \phi$.}
\item{ There exists $\psi: \R^n \longrightarrow \R_+$ convex such that
$\phi \sim \psi$.}
\end{enumerate}
\end{itemize}
\end{Prop}
{\bf Proof.} $1\im2$. Suppose $\phi$ is quasi-convex. It
suffices to show that there exists $M>0$ such that $\phi \leq M \co \phi$.
Note that the quasi-convexity of $\phi$ gives that for every $ N \geq
2$ there exists $ L_N >0$ such that for every $\{t_i\}_{i=1}^N$ in
$\R^n$ and for every $\{\lambda_i\}_{i=1}^N$ in $[0,1]$ with
$\sum_{i=1}^N \lambda_i =1$ we have:
\begin{eqnarray} \label{in:quasi}
\phi \left(
\sum_{i=1}^N \lambda_i t_i \right) \leq L_N \sum_{i=1}^N \lambda_i
\phi (t_i).
\end{eqnarray}
The proof goes by induction upon $N$. For example, for $N=3$:
\begin{eqnarray*}
\phi \left( \sum_{i=1}^3 \lambda_i t_i \right) &=& \phi \left(
\lambda_1 t_1 + (\lambda_2 + \lambda_3 ) \left(
\frac{\lambda_2}{(\lambda_2 + \lambda_3 )} t_2 +
\frac{\lambda_3}{(\lambda_2 + \lambda_3 )} t_3 \right) \right) \\
&\leq& L \left( \lambda_1 \phi(t_1) + (\lambda_2 + \lambda_3 )\phi\left(
\frac{\lambda_2}{(\lambda_2 + \lambda_3 )} t_2 +
\frac{\lambda_3}{(\lambda_2 + \lambda_3 )} t_3 \right)\right)\\
&\leq& L ( \lambda_1 \phi(t_1) + L ( \lambda_2 \phi(t_2) + \lambda_3
\phi(t_3) ))\\
&\leq& L_3 \left( \sum_{i=1}^3 \lambda_i \phi(t_i)\right) \mbox{,
where } L_3 = L^2 .
\end{eqnarray*}
Note that $\co (M\phi)= M \co\phi$. We will show $\phi \leq \co(M\phi)$
with $M=L_{n+1}$. Let $t\in\R^n$ and let $\alpha_i \geq 0 , i=1,\dots
,m$ such that $\sum_{i=1}^m \alpha_i = 1$ and $\sum_{i=1}^m \alpha_i
t_i = t$. The point $\left(t, \sum_{i=1}^m \alpha_i (L_{n+1}
\phi(t_i))\right) = \sum_{i=1}^m \alpha_i (t_i , L_{n+1} \phi(t_i))$
lies inside the convex hull of $\{(t_i , L_{n+1} \phi(t_i)) | i=1,\dots
,m\}$. Therefore, by Caratheodory's Theorem (see for example \cite{R}),
there exist $n+1$ indices
$i_1 , \dots , i_{n+1}$ and $ \lambda_1 , \dots ,\lambda_{n+1}\geq 0$
with $\sum_{i=1}^{n+1} \lambda_i =1$ such that:
\begin{eqnarray} \label{Cara}
& &t=\sum_{j=1}^{n+1} \lambda_j t_{i_j} \mbox{ and } \nonumber\\
& &\sum_{j=1}^{n+1} \lambda_j L_{n+1} \phi (t_{i_j}) \leq \sum_{i=1}^m
\alpha_i L_{n+1} \phi (t_i ).
\end{eqnarray}
By applying (\ref{in:quasi}) for $N=n+1$ we see that
\[ \phi(t) \leq L_{n+1} (\sum_{j=1}^{n+1} \lambda_j \phi (t_{i_j})). \]
Hence, by (\ref{Cara}) we get
\[ \phi(t) \leq \sum_{i=1}^m \alpha_i L_{n+1} \phi (t_i ).\]
By taking the infimum over all convex combinations $t=\sum_{i}
\alpha_i t_i$ we get $\phi(t) \leq \co(L_{n+1}\phi)(t)$. QED.

\noindent
$2\im3$ is trivial, just let $\psi=\co\phi$.

\noindent
$3\im1$. Suppose $\psi$ is convex and let $M>0$ such that
$\frac{1}{M}\psi(x) \leq\phi(x) \leq M \psi(x)$ for all $x\in
\R^n$. Then:
\begin{eqnarray*}
\phi(\lambda_1 x_1 + \lambda_2 x_2 ) &\leq& M \psi(\lambda_1 x_1 +
\lambda_2 x_2 )
\leq  M(\lambda_1 \psi(x_1 ) +\lambda_2 \psi(x_2 ))\\
& & \leq M^2
(\lambda_1 \phi(x_1 ) + \lambda_2 \phi(x_2 )).
\end{eqnarray*}
Thus $\phi$ is quasi-convex. \hfill $\Box$

\smallskip
\noindent
The next theorem will give examples of quasi-convex functions on
$\R^2$ (which are not convex). These examples will play an important
role in the next section. We recall that $\phi$ is an Orlicz function
if it is a convex, non-decreasing function on $[0,\infty)$ such that
$\phi(0)=0$ and $\lim_{t\rightarrow\infty} \phi(t)=\infty$. For more
information on Orlicz spaces see \cite{LT1}. The functions we shall
consider will be {\it finite valued} and {\it non-degenerate}, that is
$0$ only at $0$. We say $\phi$ satisfies the {\it $\Delta_2$ condition at
zero } if $\limsup_{x\rightarrow 0}\frac{\phi(2x)}{\phi(x)}<\infty$ and,
respectively, $\phi$ satisfies the {\it $\Delta_2$ condition } if there
exists $C>0$ such that for all $x \geq 0 \mbox{, }\phi(2x)\leq
C\phi(x)$. We will extend an Orlicz function on the whole real
line by $\phi(x)=\phi(-x)$ if $x<0$ and, abusing the language, we will
still call the extension an Orlicz function. We start with the
following simple
\begin{Obs} \label{O:cond}
Let $\phi$ be an Orlicz function such that
\begin{equation} \label{E:typep}
\exists p>1 \mbox{, } \exists M>0 \mbox{ such that } \forall \lambda
\in (0,1]\mbox{, } \forall s>0 \mbox{, } \frac{\phi(\lambda
s)}{\lambda^p \phi(s)}\leq M .
\end{equation}
Then
\begin{equation} \label{E:type}
\exists M'>0 \mbox{ such that } \forall \lambda \in (0,1],\forall
y>0 \mbox{ we have }
\frac{\phi(\lambda |\log(\lambda)|y)}{\lambda\phi(y)}\leq M' .
\end{equation}
\end{Obs}
Indeed, suppose (\ref{E:typep}) holds. Note that for $\lambda \in
[0,1]$, $\lambda |\log\lambda|\in
[0,\frac{1}{e}]$. Therefore if $\lambda\in (0,1]$ and $y>0$ we have:
\begin{eqnarray*}
\frac{\phi(\lambda|\log\lambda|y)}{\lambda\phi(y)}=
\frac{\phi(\lambda|\log\lambda|y)}{\lambda^p |\log\lambda|^p \phi(y)}
\cdot \lambda^{p-1}|\log\lambda|^p \leq M S < \infty
\end{eqnarray*}
where $S= \sup_{\lambda \in [0,1]} \lambda^{p-1}|\log\lambda|^p$.
We are now ready to state the main result  of this section.
\begin{Thm} \label{T:ex}
Let $\phi$ be an Orlicz function satisfying {\rm (\ref{E:typep})}
and the $\Delta_2$ condition. Let $\theta : \R \rightarrow \R$ be a
Lipschitz map. Then
$\Phi : \R^2 \rightarrow \R_+$ defined by
\[ \Phi(x,y)=\left\{ \begin{array}{ll}
			\phi(y) + \phi (x- y
			\theta(\log\frac{1}{|y|})) &\mbox{, if } y \neq
			0 \\
			\phi(x) &\mbox{, if } y=0
		    \end{array}
			\right. \]
is quasi-convex.
\end{Thm}
{\bf Proof.} By Observation \ref{O:cond} the hypothesis
(\ref{E:typep}) gives (\ref{E:type}).
Using the $\Delta_2$ condition and the increasingness
of $\phi$ one can easily prove that there exists $C>0$ such that for
all $x,y\in\R$  we have
\begin{equation} \label{E:indel}
\phi(x+y) \leq C(\phi(x)+\phi(y))
\end{equation}
and that for all $B>0$ there
exists $C_B >0$ such that for all $x \geq 0$
\begin{equation} \label{E:indelb}
\phi(Bx) \leq C_B \phi(x).
\end{equation}
Let $t_i =(x_i ,  y_i ) \in \R^2$ and $\lambda_i \in [0,1] \mbox{, }
i=1,2 \mbox{ with }\lambda_1 + \lambda_2 =1$. Without loss of
generality we may assume that $\lambda_1 y_1 \neq 0$ and $\lambda_2
y_2 \neq 0$. Then
\begin{eqnarray*}
\Phi\left(\sum_{i=1}^2 \lambda_i t_i\right)&=&
\phi\left(\sum_{i=1}^2 \lambda_i y_i\right)+
\phi\left(\sum_{i=1}^2 \lambda_i x_i - \sum_{i=1}^2 \lambda_i y_i \theta
\left( \log\frac{1}{|\sum_{i=1}^2 \lambda_i y_i |}\right)\right) \\
&=&\phi\left(\sum_{i=1}^2 \lambda_i y_i\right)+\phi\left(\sum_{i=1}^2
\lambda_i x_i
-\sum_{i=1}^2 \lambda_i y_i \theta \left( \log\frac{1}{|\lambda_i y_i
|}\right)\right.+ \\
& &\sum_{i=1}^2 \lambda_i y_i \theta \left( \log\frac{1}{|\lambda_i y_i
|}\right)
\left. -  \sum_{i=1}^2 \lambda_i y_i \theta
\left( \log\frac{1}{|\sum_{i=1}^2 \lambda_i y_i |}\right)\right) \\
&\leq& \phi\left(\sum_{i=1}^2 \lambda_i y_i\right) + C\phi\left(\sum_{i=1}^2
\lambda_i x_i -\sum_{i=1}^2 \lambda_i y_i \theta
\left( \log\frac{1}{| \lambda_i y_i |}\right)\right) \\
& &+C\phi\left(\sum_{i=1}^2 \lambda_i y_i \theta \left(
\log\frac{1}{|\lambda_i y_i
|}\right)- \sum_{i=1}^2 \lambda_i y_i \theta
\left( \log\frac{1}{|\sum_{i=1}^2 \lambda_i y_i |}\right)\right) \\
\end{eqnarray*}
using (\ref{E:indel}). For the last term in the sum we apply the inequality
\begin{eqnarray*}
|t\theta(\log\frac{1}{|t|}) +
s\theta(\log\frac{1}{|s|})-(t+s)\theta(\log\frac{1}{|t+s|})|\leq K
(|s|+|t|)
\end{eqnarray*}
where $K$ is the Lipschitz constant of $\theta$. This inequality shows
that the map $t \mapsto t\theta(\log\frac{1}{|t|})$ is quasi-additive
(see \cite{KP}, Theorem 3.7). Since $\phi$ is increasing on the
positive axis, we obtain:
\begin{eqnarray*}
& &\Phi(\sum_{i=1}^2 \lambda_i t_i)\leq\phi(\sum_{i=1}^2 \lambda_i
y_i) + C\phi\left(\sum_{i=1}^2
\lambda_i x_i -\sum_{i=1}^2  \lambda_i y_i \theta
\left( \log\frac{1}{| \lambda_i y_i |}\right)\right)\\
& &+C\phi(K\sum_{i=1}^2 \lambda_i |y_i |)\\
&\leq&(1+C C_K )\sum_{i=1}^2 \lambda_i \phi(y_i )+
 C\phi\left(\sum_{i=1}^2
\lambda_i x_i -\sum_{i=1}^2  \lambda_i y_i \theta
\left( \log\frac{1}{| \lambda_i y_i |}\right)\right)
\end{eqnarray*}
by using the convexity of $\phi$ and (\ref{E:indelb}). Thus
\begin{eqnarray*}
& &\Phi\left(\sum_{i=1}^2 \lambda_i t_i\right)\leq(1+C C_K )\sum_{i=1}^2
\lambda_i \phi(y_i )+ C\phi\left(\sum_{i=1}^2
\lambda_i x_i - \right.\\
& &\left. - \sum_{i=1}^2 \lambda_i y_i \theta
\left( \log\frac{1}{|y_i |}\right)
+\sum_{i=1}^2 \lambda_i y_i \theta
\left( \log\frac{1}{|y_i |}\right)-\sum_{i=1}^2  \lambda_i y_i \theta
\left( \log\frac{1}{| \lambda_i y_i |}\right)\right)\\
&\leq&(1+C C_K )\sum_{i=1}^2 \lambda_i \phi(y_i ) + C^2
\phi\left(\sum_{i=1}^2
\lambda_i x_i - \sum_{i=1}^2 \lambda_i y_i \theta
\left( \log\frac{1}{|y_i |}\right)\right)\\
& & + C^2 \phi\left(\sum_{i=1}^2 \lambda_i y_i \theta
\left( \log\frac{1}{|y_i |}\right)-\sum_{i=1}^2  \lambda_i y_i \theta
\left( \log\frac{1}{| \lambda_i y_i |}\right)\right)\\
&\leq&(1+C C_K )\sum_{i=1}^2 \lambda_i \phi(y_i )+
C^2\sum_{i=1}^2\lambda_i\phi\left(x_i -y_i
\theta\left(\log\frac{1}{|y_i|}\right)\right) \\
& & +C^2\phi\left(\sum_{i=1}^2 \lambda_i |y_i|
\left|\theta\left(\log\frac{1}{|y_i|}\right)-\theta\left(\log\frac{1}{|\lambda_i
y_i|}\right)\right|\right)\\
&\leq&(1+C C_K )\sum_{i=1}^2 \lambda_i \phi(y_i )+
C^2\sum_{i=1}^2\lambda_i\phi\left(x_i -y_i
\theta\left(\log\frac{1}{|y_i|}\right)\right) \\
& & +C^3\sum_{i=1}^2\phi\left( \lambda_i|y_i| K
\left|\log\frac{1}{|y_i|}-\log\frac{1}{|\lambda_i
y_i|}\right|\right)\\
&\leq&(1+C C_K )\sum_{i=1}^2 \lambda_i \phi(y_i )+
C^2\sum_{i=1}^2\lambda_i\phi\left(x_i -y_i
\theta\left(\log\frac{1}{|y_i|}\right)\right) \\
& & +C^3 C_K \sum_{i=1}^2 \phi( \lambda_i |y_i||\log\lambda_i|).
\end{eqnarray*}
By applying (\ref{E:type}) to the last term we obtain
\begin{eqnarray*}
\Phi\left(\sum_{i=1}^2 \lambda_i t_i\right)&\leq&(1+C C_K )\sum_{i=1}^2
\lambda_i \phi(y_i )+
C^2\sum_{i=1}^2\lambda_i\phi\left(x_i -y_i
\theta\left(\log\frac{1}{|y_i|}\right)\right) \\
& &+C^3 C_K M' \sum_{i=1}^2 \lambda_i\phi(|y_i|)\\
&\leq& \max(1+C C_K +C^3 C_K M' , C^2)\sum_{i=1}^2 \lambda_i\Phi(t_i)
\end{eqnarray*}
which ends the proof. \hfill $\Box$
\smallskip
\noindent
\begin{Rmk}
If $\ell_\phi$ is an Orlicz space with type greater than $1$ then
there exists an Orlicz function $\tilde{\phi}$ satisfying
{\rm(\ref{E:typep})} and the $\Delta_2$ condition such that $\tilde{\phi}$
coincides with $\phi$ on $[0,1]$.
\end{Rmk}
Indeed, it is well-known that the space $\ell_\phi$ has non-trivial
type if and only if $\alpha_\phi >1$ and $\beta_\phi<\infty$, where
$\alpha_\phi$ and $\beta_\phi$ are the lower and the upper indices:
\begin{eqnarray*}
& & \alpha_\phi = \sup\{q; \sup_{0< \lambda, t \leq 1}
\frac{\phi(\lambda t)}{\phi(\lambda)t^q} <\infty \} \mbox{ and }\\
& & \beta_\phi = \inf\{q; \inf_{0< \lambda, t \leq 1}
\frac{\phi(\lambda t)}{\phi(\lambda)t^q}>0\}
\end{eqnarray*}
(cf. \cite{LT2} p.140 and \cite{LT1} p.143). Moreover, $\beta_\phi<\infty$ is
equivalent to
$\phi$ satisfying the $\Delta_2$ condition at zero (see \cite{LT1}).
Note that  $\alpha_\phi >1$ means:
\begin{equation} \label{eqtp}
\exists p >1 \mbox{, } \exists M>0 \mbox{ such that }\forall \lambda
\in (0,1]\mbox{, } \forall s\in(0,1] \mbox{, } \frac{\phi(\lambda
s)}{\lambda^p \phi(s)}\leq M .
\end{equation}
Define
\[\tilde{\phi}(x)=\left\{\begin{array}{ll}
				\phi(x) &\mbox{, if } x\leq 1 \\
				\phi(1)x^q &\mbox{, if } x>1
			\end{array}
		     	\right. \]

\noindent
where $q=\max\{\frac{\phi'(1)}{\phi(1)},p\}$ and $\phi'(1)$ denotes the
left derivative of $\phi$ at $1$. Clearly $\tilde{\phi}$ is
Orlicz. Since $\phi$ satisfies the $\Delta_2$ condition at zero,
$\tilde{\phi}$ satisfies the $\Delta_2$ condition (note that we don't
use the full assumption of the existence of type for this part of the
argument).  Let us check that
$\tilde{\phi}$ satisfies (\ref{E:typep}). Let $\lambda\in(0,1)$ and
$s\in(0,\infty)$. If $s\leq1$ the inequality in (\ref{E:typep}) is
given by (\ref{eqtp}). If $s>1$ we consider two cases: if $\lambda s
\leq 1$ then
\begin{eqnarray*}
\frac{\tilde{\phi}(\lambda s)}{\lambda^p \tilde{\phi}(s)}=
\frac{\phi(\lambda s)}{\lambda^p \phi(1) s^q} \leq \frac{(\lambda s)^p M
\phi(1)}{\lambda^p \phi(1) s^q}=\frac{M}{s^{q-p}}\leq M
\end{eqnarray*}
and if $\lambda s> 1$ then
\begin{eqnarray*}
\frac{\tilde{\phi}(\lambda s)}{\lambda^p \tilde{\phi}(s)}=
\frac{\phi(1) (\lambda s)^q}{\lambda^p \phi(1)s^q}= \lambda^{q-p}\leq
1.
\end{eqnarray*}

\noindent
Clearly if an
Orlicz function $\phi$ satisfies (\ref{E:typep}) and the $\Delta_2$
condition, then $\ell_\phi$ has non-trivial type. Finally, we note
that Theorem \ref{T:ex} implies that for an Orlicz space  $\ell_\phi$
with non-trivial type, there exists an Orlicz function $\tilde{\phi}$
generating  $\ell_\phi$
such that $\Phi : \R^2 \rightarrow \R_+$ defined by
\[ \Phi(x,y)=\left\{ \begin{array}{ll}
			\tilde{\phi}(y) + \tilde{\phi} (x- y
			\theta(\log\frac{1}{|y|})) &\mbox{, if } y \neq
			0 \\
			\tilde{\phi}(x) &\mbox{, if } y=0
		    \end{array}
			\right. \]
is quasi-convex.

\section{Twisted Sums and Fenchel-Orlicz Spaces} \label{S:TSFO}
\noindent
Let $c_{00}$ denote the space of real sequences with finite support. Let
$\theta : \R \rightarrow \R$ be a
Lipschitz map. Let $\phi$
be an Orlicz function satisfying the $\Delta_2$ condition such that
$\ell_\phi$ has non-trivial type. Let $\|\cdot\|_{\phi}$ denote the
norm of the Orlicz space $\ell_\phi$:
\[\|(x_n)_n \|_\phi = \inf \{ \rho > 0 : \sum_{n}
\phi\left(\frac{x_n}{\rho}\right)\leq 1 \}. \]
Let $F: c_{00} \longrightarrow c_{00}$
be defined by:
\[ (F(y_m)_m)_n= \left\{ \begin{array}{ll}
				 y_n \theta \left(
				 \log
				 \frac{\|(y_m)_m\|_\phi}{|y_n|}\right)
				 &\mbox{, if } y_n \neq 0 \\
				 0 &\mbox{, if } y_n = 0
			\end{array}
				\right. \]
It is proved in \cite{KP} that F is a quasi-linear map, i.e. for all
$ \lambda  \in \R$ and for all $x,y\in c_{00}$ we have:
\begin{eqnarray*}
& & F(\lambda y)=\lambda F(y)  \mbox{ and}\\
& & \|F(x+y)-F(x)-F(y)\| \leq c (\|x\|+\|y\|)
\end{eqnarray*}
where $c$ is a constant independent of $x$ and $y$.
We define a quasi-norm  on $c_{00}\times c_{00}$ by
\[\|(x_n, y_n)_n\| = \|(y_n)_n\|_\phi +\|(x_n)_n - F((y_m)_m)_n \|_\phi
.\]
The twisted sum $\ell_\phi \bigoplus_F \ell_\phi$ is defined as the
completion of $c_{00}\times c_{00}$
with respect to the quasi-norm
$\|\cdot\|$. In other words, $\ell_\phi \bigoplus_F \ell_\phi$
consists of all sequences
$(x_n,y_n)_n$ such that $\|(x_n, y_n)_n\|<\infty$. The fact that
$\ell_\phi \bigoplus_F \ell_\phi$ is a Banach space follows from
Theorem 2.6 in \cite{K1} which implies that a twisted sum of two
B-convex Banach spaces is (after renorming) a B-convex Banach space
and Pisier's result \cite{P} that a Banach space $X$ has type greater than $1$
if and only if it is B-convex.
\bigskip
\noindent
\begin{Def}[\cite{T}] \label{D:Yfn}
 A Young's function on $\R^n$ is
an even, convex function $\Phi : \R^n \rightarrow \overline{\R}_+$
with $\Phi(0)=0$ and $\lim_{t\rightarrow\infty}\Phi(tx)=\infty$
for all $x\in\R^n\setminus\{0\}$.
\end{Def}
We set (cf. \cite{LT1})
\begin{eqnarray} \label{D:ell}
\ell_\Phi = \{ (x_k ^1,\dots ,x_k^n)_k :\exists \rho >0 \mbox{ such that }
\sum_{k}\Phi\left(\frac{1}{\rho}(x_k ^1,\dots ,x_k^n)\right) <\infty\}
\end{eqnarray}
and for $(x_k ^1,\dots ,x_k^n)_k \in \ell_\Phi$ we define
\begin{eqnarray*}
\|(x_k ^1,\dots ,x_k^n)_k\|_\Phi = \inf\{\rho>0 :
\sum_{k}\Phi\left(\frac{1}{\rho}(x_k ^1,\dots ,x_k^n)\right)  \leq 1\}.
\end{eqnarray*}
Then $\ell_\Phi$ is a vector space and $(\ell_\Phi, \|\cdot\|_\Phi)$ is called
a Fenchel-Orlicz space. If $\Phi$ is {\it finite} on $\R^n$, $\ell_\Phi$ is
complete in $\| \cdot \|$ (see Corollary 2.23 in \cite{T}).
For a detailed study of (more general) Fenchel-Orlicz spaces and their
completeness we refer to Turett \cite{T}. Note that for $n=1$ we
retrieve the Orlicz spaces.
We also define $h_\Phi$ to be the vector subspace of $\ell_\Phi$ consisting
of all sequences $(x_k^1,\dots , x_k^n)_k$ such that
$\sum_{k}\Phi(\frac{1}{\rho}(x_k^1,\dots, x_k^n))<\infty$
for every $\rho>0$.
With abuse of notation we will use (\ref{D:ell}) to define $\ell_\Phi$
for any quasi-convex function $\Phi : \R^n \rightarrow \R_+$; similarly for
$h_\Phi$. We will say that a quasi-convex map $\Phi : \R^n \rightarrow
\R_+$ satisfies the {\it $\Delta_2$ condition} if there exists $M>0$
such that $\Phi(2x)\leq M \Phi(x)$ for all $x\in \R^n$. Note that if
$\Phi : \R^n \rightarrow \R_+$ is a quasi-convex even function with
$\Phi(0)=0$ and $\lim_{t\rightarrow\infty}\Phi(tx)=\infty$
for all $x\in\R^n\setminus\{0\}$ then Proposition \ref{P:quasi}
implies that, as sets,
\begin{eqnarray} \label{E:ega}
\ell_\Phi = \ell_{\co\Phi} \mbox{ and } h_\Phi=h_{\co\Phi}.
\end{eqnarray}

\bigskip
\noindent
The main result of this section is:
\begin{Thm} \label{T:twistOFO}
If $\ell_\phi$ is an Orlicz space with non-trivial type then the twisted sum
$\ell_\phi \bigoplus_F \ell_\phi$ is a Fenchel-Orlicz space on
$\R^2$. More precisely,there exists a Young's function $\Psi$ on
$\R^2$ such that $\ell_\phi \bigoplus_F \ell_\phi = \ell_\Psi$
(as sets) and the identity map is an isomorphism.
\end{Thm}

\noindent
The rest of this section will be devoted to the proof of this result.
By Theorem \ref{T:ex} and the remarks following it we see that, without
loss of generality, we may assume that $\phi$ satisfies the $\Delta_2$
condition and that the map $\Phi : \R^2\longrightarrow \R_+$ defined by
\[ \Phi(x,y)=\left\{ \begin{array}{ll}
			\phi(y) + \phi (x- y
			\theta(\log\frac{1}{|y|})) &\mbox{, if } y \neq
			0 \\
			\phi(x) &\mbox{, if } y=0
		    \end{array}
			\right. \]
is quasi-convex. We shall show that the function $\Psi=\co\Phi$ is a
Young's function on $\R^2$ with the property mentioned in the theorem.

\bigskip
\noindent
We first prove the set equality between the two spaces. We start with
\begin{Rmk}
 $\Phi$ satisfies the $\Delta_2$ condition.
\end{Rmk}
Indeed, for $y\neq0$ we have:
\begin{eqnarray*}
& &\phi\left(2x-2y\theta\left(\log\frac{1}{2|y|}\right)\right) \leq C
\phi\left(x-y\theta\left(\log\frac{1}{2|y|}\right)\right)\\
&=&C\phi\left(x-y\theta\left(\log\frac{1}{|y|}\right)+
y\theta\left(\log\frac{1}{|y|}\right)-
y\theta\left(\log\frac{1}{2|y|}\right)\right)\\
&\leq& C^2\phi\left(x-y\theta\left(\log\frac{1}{|y|}\right)\right)+
C^2\phi(yK\log2)\\
&\leq& C^2\phi\left(x-y\theta\left(\log\frac{1}{|y|}\right)\right)+
C^2 C_{K\log2}\phi(y)
\end{eqnarray*}
where $K$ is the Lipschitz constant of $\theta$ while $C$ and
$C_{K\log2}$ are given by (\ref{E:indel}) and (\ref{E:indelb})
respectively. Note that if a quasi-convex function $\psi:\R^n
\rightarrow \R$ satisfies the $\Delta_2$ condition then
$\ell_\psi=h_\psi$ (cf.\cite{LT1} Proposition 4.a.4). Therefore:
\begin{equation} \label{E:lfi=hfi}
\ell_\Phi = h_\Phi
\end{equation}

\bigskip
\noindent
The following notation will simplify further computations. For a given
sequence $(x_j ,y_j )_j$ let
\begin{eqnarray*}
S(k) =  \sum_{j} \phi(y_j) + \sum_{j} \phi\left(x_j - y_j\theta\left(
\log{\frac{k}{|y_j|}}\right)\right)
\end{eqnarray*}
for $k>0$. It is easy to see that
\[ \ell_\Phi = \{ (x_j, y_j )_j |\mbox{ there exists }\rho >0 ,
S( \rho )< \infty\}\]
and
\[  h_\Phi = \{ (x_j, y_j )_j | \mbox{ for all } \rho >0 , S( \rho )< \infty
\}.\]
Indeed, if $(x_j, y_j )_j \in  \ell_\Phi$ there exists $\rho>0$ such
that
\[\sum_{j} \phi\left(\frac{y_j}{\rho}\right) + \sum_{j}
\phi\left(\frac{x_j -
y_j\theta\left(\log{\frac{\rho}{|y_j|}}\right)}{\rho}\right) <
\infty. \]
But then, since $\phi$ satisfies the $\Delta_2$ condition,
\begin{eqnarray*}
S(\rho)&=&\sum_{j} \phi\left(\rho \frac{y_j}{\rho}\right) +
\sum_{j} \phi\left(\rho \frac{x_j -
y_j\theta\left(\log{\frac{\rho}{|y_j|}}\right)}{\rho}\right) \\
&\leq& C_{\rho} \left(\sum_{j} \phi\left(\frac{y_j}{\rho}\right) + \sum_{j}
\phi\left(\frac{x_j -
y_j\theta\left(\log{\frac{\rho}{|y_j|}}\right)}{\rho}\right)\right) < \infty.
\end{eqnarray*}
Conversely, if $(x_j, y_j )_j$ is such that $ S( \rho )< \infty$ then
\[ \sum_{j} \phi\left(\frac{y_j}{\rho}\right) + \sum_{j}
\phi\left(\frac{x_j -
y_j\theta\left(\log{\frac{\rho}{|y_j|}}\right)}{\rho}\right) \leq
C_{\frac{1}{\rho}} S(\rho) < \infty \]
and thus  $(x_j, y_j )_j \in  \ell_\Phi$.

\smallskip
\noindent
Moreover, note that for $\|(y_j)_j\|_\phi >0$ we have:
\begin{equation} \label{E:equi}
 \|(x_j,y_j)_j\| < \infty \mbox{ if and only if (}\|(y_j)_j\|_\phi< \infty
\mbox{ and) }
S(\|(y_j)_j\|_\phi)< \infty .
\end{equation}

\smallskip
\noindent
We now show that $\ell_\phi \bigoplus_F \ell_\phi = \ell_\Psi$ as sets.
Let $(x_j, y_j)_j \in\ell_\phi \bigoplus_F \ell_\phi .$  Then $\|(x_j
,y_j )\|<\infty$, which implies $(y_j)_j \in \ell_\phi$. If
$\|(y_j)_j\|_\phi > 0$ then by (\ref{E:equi}) we get
$S(\|(y_j)_j\|_\phi)<\infty$. This shows
that  $(x_j, y_j)_j \in\ell_\Phi$. If $\|(y_j)_j\|_\phi =0 $ then
$(x_j)_j\in \ell_\phi$ and again  $(x_j, y_j)_j \in\ell_\Phi$. Hence, by
(\ref{E:ega}), $(x_j, y_j)_j \in\ell_\Psi$.  Conversely if  $(x_j,
y_j)_j \in\ell_\Psi$, by (\ref{E:ega}) and (\ref{E:lfi=hfi}),
$(x_j,y_j)_j \in h_\Phi$, which implies $(y_j)_j \in \ell_\phi$. If
$(y_j)_j \neq 0$ then $S(\|(y_j)_j\|_\phi)<\infty$. Therefore $ \|(x_j
,y_j)\| <\infty$ and  $(x_j, y_j)_j \in\ell_\phi \bigoplus_F \ell_\phi
.$ If $(y_j)_j = 0$ then $(x_j)_j \in \ell_\phi$ and again $(x_j,
y_j)_j \in\ell_\phi \bigoplus_F \ell_\phi.$

\smallskip
\noindent
Note that $\Psi$ is a finite Young's function and thus $\ell_\Psi$ is
a Banach space. Indeed, we only need to show that
\[ \lim_{t \rightarrow \infty} \Phi (t(x,y)) = \infty \mbox{, for all }(x,y)
\in \R^2
\setminus \{ (0,0) \} \]
since then Proposition \ref{P:quasi} will give the same result for $\Psi$.
Let $(x,y)\neq(0,0)$. If $y\neq 0$ then $\Phi(t(x,y))\geq
\phi(ty) \rightarrow \infty$ as $ t \rightarrow \infty$ since $\phi$
is an Orlicz function. If $y=0$ then $x\neq0$ and
$\Phi(t(x,y))=\phi(tx) \rightarrow \infty$ as $ t \rightarrow \infty$
since $\phi$ is an Orlicz function.

\smallskip
\noindent
The next two propositions will show that the identity mapping
is an isomorphism between $\ell_\phi \bigoplus_F \ell_\phi$ and
$\ell_\Psi$.
\begin{Prop}
Let $X$ be a sequence space, complete in $\|\cdot\|_1$ and $\|\cdot\|_2$,
such that the coordinate functionals are continuous. Then the identity
$\imath : (X,\|\cdot\|_1) \rightarrow (X,\|\cdot\|_2)$ is an isomorphism.
\end{Prop}
{\bf Proof.} It is easy to see that $(X,\|\cdot\|_1 +\|\cdot\|_2)$ is
complete. Therefore the identity maps
\[ \imath_1 :(X,\|\cdot\|_1+\|\cdot\|_2) \rightarrow
(X,\|\cdot\|_1)\mbox{ and }
\imath_2 :(X,\|\cdot\|_1+\|\cdot\|_2) \rightarrow (X,\|\cdot\|_2) \]
are continuous and hence, by the Inverse Mapping Theorem,
isomorphisms. Therefore $\imath=\imath_2 \circ \imath_1 ^{-1}$ is an
isomorphism. \hfill $\Box$

\bigskip
\noindent
\begin{Prop}
The coordinate functionals on $\ell_\phi \bigoplus_F \ell_\phi$ and
$\ell_\Psi$ are continuous.
\end{Prop}
{\bf Proof.} In both cases we will show that projections $P_i((x_j
 ,y_j )_j) =(x_i , y_i)$ are continuous for all $i$. The result will follow
immediately since the coordinate functionals on 2-dimensional Banach
spaces are continuous.

For $\ell_\Psi$ note that if $\|(x_j ,y_j )_j\|_\Psi =1$ then by the
continuity of $\Psi$
we get $\sum_{j} \Psi (x_j ,y_j) \leq 1$ and hence $\Psi(x_i ,y_i)
\leq 1$ for all $i$. Therefore, for all $i$,
$(x_i,y_i)\in\{(x,y)\in\R^2| \Psi(x,y) \leq 1\}$ which is a bounded
set and thus the projection $P_i$ is continuous since all
norm-topologies on a 2-dimensional Banach space are equivalent.

For $\ell_\phi \bigoplus_F \ell_\phi$ we show that there exists $M$
such that if $\|(x_j ,y_j )_j\| \leq 1$ then $\Phi(x_i ,y_i)\leq
M$. The boundedness of the set $\{(x,y)\in\R^2 | \Phi(x,y) \leq M\}$
finishes the proof as before. Indeed, note that if $\|(x_j ,y_j )_j\|
\leq 1$ then $S(\|(y_j)_j\|_\phi) \leq1$. Moreover, if $C$ is given by
(\ref{E:indel}) and $K$ is the Lipschitz constant of $\theta$ then
\begin{eqnarray*} \label{ineq}
\Phi(x_i,y_i)&\leq&S(1)=\sum_j \phi(y_j) + \sum_j \phi\left(x_j -y_j
\theta\left(\log\frac{\|(y_n)_n\|_\phi}{|y_j|}\right)\right.+\\
& &\left.+y_j\theta\left(\log\frac{\|(y_n)_n\|_\phi}{|y_j|}\right)-
 y_j\theta\left(\log\frac{1}{|y_j|}\right)\right) \\
&\leq&\sum_j \phi(y_j)+C\sum_j \phi\left(x_j -y_j
\theta\left(\log\frac{\|(y_n)_n\|_\phi}{|y_j|}\right)\right)+\\
& & +C\sum_j\phi\left(|y_j K\log\|(y_n)_n\|_\phi |\right)\\
&\leq&(1+C)S(\|(y_n)_n\|_\phi)+C\sum_j\phi\left(y_j
K\log\|(y_n)_n\|_\phi \right)\\
&\leq& (1+C) +C \sum_j \phi\left(\frac{y_j}{\|(y_n)_n\|_\phi}\|(y_n)_n\|_\phi
 K \log\|(y_n)_n\|_\phi\right) \\
&\leq& 1+C + M'C \sum_j
\phi\left(\frac{y_j}{\|(y_n)_n\|_\phi}K\right)\|(y_n)_n\|_\phi
\end{eqnarray*}
where $M'$ is given by (\ref{E:type}). Thus, if $C_K$ is given by
(\ref{E:indelb}), we have:
\begin{eqnarray*}
\Phi(x_i,y_i)&\leq& 1+C + M'C C_K \sum_j
\phi\left(\frac{y_j}{\|(y_n)_n\|_\phi}\right)\leq 1+C + M'CC_K.
\end{eqnarray*}
The proof of the proposition is complete. \hfill $\Box$

\noindent
This concludes the proof of Theorem \ref{T:twistOFO}.

\smallskip
\noindent
In particular, by choosing $\phi(x)= |x|^p$ for $ 1<p<\infty$ and
$\theta$ to be the identity map, we see that the spaces $Z_p$
introduced in \cite{KP} can be viewed as Fenchel-Orlicz
spaces. We end this section with two questions which arise naturally,
in view of Theorem \ref{T:twistOFO}:
For what Banach spaces $X$ can a twisted sum of $X$ with itself be
represented as a Fenchel-Orlicz space? Note that this can not be done
for $X=\ell_1$ as $\ell_1 \bigoplus_F \ell_1$ is not a Banach space.
If $\ell_\phi$ is an Orlicz space with non-trivial type for which
quasi-linear maps $G$ is the
twisted sum  $\ell_\phi \bigoplus_G \ell_\phi$ a Fenchel-Orlicz space?

\bigskip
\noindent
\section{Fenchel-Orlicz spaces with property (M)} \label{S:FOM}
Recall the definition of property (M) \cite{K4} (see also \cite{HWW}):
\begin{Def}
A Banach space $X$ has property (M) if whenever $u,v \in X$ with $\|u
\| = \|v\|$
and $(x_n)_n$ is a weakly null sequence in $X$ then
\[\limsup_{n \rightarrow \infty} \| u + x_n \| =\limsup_{n \rightarrow
\infty} \| v + x_n \| \]
\end{Def}
A large class of spaces with property (M) can be generated as follows:
Let $(n_k )_k$ be a sequence of natural numbers.
For every $k$ let $N_k$ be a norm on $ \R^{n_k + 1}$ such that
\[ 0\leq x_0 \leq x_0 ' \Rightarrow  N_k ( x_0 , x_1 ,\dots ,
x_{n_k}) \leq
N_k ( x_0 ', x_1 ,\dots , x_{n_k}) \]
and
\[N_k (1,0,\dots,0)=1. \]
Define inductively a sequence of norms on $ \R ^ {\sum_{i=1}^k n_i }$ by:
\[ N_1 \ast N_2 (x_1 , x_2 , \dots , x_{n_1 + n_2} ) = N_2 ( N_1 ( 0 ,
x_1 , \dots , x_{n_1}), x_{n_1 +1}, \dots ,
x_{n_1 + n_2}) \]
and once $N_1 \ast \dots \ast N_{k-1}$ is defined,
\begin{eqnarray*}
& & N_1 \ast \dots \ast N_k ( x_1 , \dots , x_{ \sum_{i=i}^k n_i }) = \\
& & N_k ( N_1 \ast \dots \ast N_{k-1} ( x_1 , \dots, x_{
\sum_{i=i}^{k-1} n_i} ), x_{ \sum_{i=i}^{k-1} n_i +1},\dots , x_{
\sum_{i=i}^k n_i })
 \end{eqnarray*}
It can be easily checked that each $N_1 * \dots * N_k $ is a norm.
For a sequence of finite sequences $\xi = ((\xi_i )_{i=1}^{n_1},(\xi_i
)_{i=n_1 + 1}^{n_2}, \dots ,(\xi_i )_{i=n_k +1}^{n_{k+1}}, \dots )$ let
\[ \| \xi \|_{\tilde{\Lambda}(N_k)} = \sup_{k} (N_1 \ast \dots \ast
N_k)(\xi_1, \dots , \xi_{\sum_{i=1}^k n_i}) \]
and let $\tilde{\Lambda}(N_k)$ be the space of all sequences of finite
sequences $\xi$ such that $\| \xi \|_{\tilde{\Lambda}(N_k)} <
\infty$. Then $\| \cdot \|_{\tilde{\Lambda}(N_k)}$ is a norm and
$\tilde{\Lambda}(N_k)$ is a Banach space. Define $\Lambda (N_k)$ to be
the closed linear span of the basis vectors $(e_k )_k$ in
$\tilde{\Lambda}(N_k)$. A simple gliding hump argument shows that
$\Lambda (N_k)$ has property (M) (see \cite{K4}). The above technique
is used in \cite{K4} to show that the closed linear span of the basis
of modular spaces can be renormed to have property (M). If $N_k =N$
for all $k$ we write $\tilde{\Lambda}(N)$ for $\tilde{\Lambda}(N_k)$
and $\Lambda (N)$ for $\Lambda (N_k)$.

\smallskip
\noindent
For the rest of the section $n \in \N$ will be fixed and for
Fenchel-Orlicz spaces $\ell_\Phi$ on $\R^n$ we shall assume that
{\it the Young's function $\Phi$ is finite and $0$ only at $0$}. Our
main result in this section is the following
\begin{Thm} \label{T:FOproM}
Every Fenchel-Orlicz space $h_\Phi$ on $\R^n$ can be equivalently\\
renormed to have property (M).
\end{Thm}
The theorem will be proved once we show that if $\Phi : \R^n
\rightarrow \R_+$ is a Young's function there exists a norm $N$ on
$\R^{n+1}$ such that $\ell_\Phi =\tilde{\Lambda} (N)$ (and thus
$h_\Phi =\Lambda (N)$ ). Sufficient conditions for
this last claim are given in the following
\begin{Lem} \label{L:suff}
If $\Phi$ is a Young's function on $\R^n$ and N is a norm on
$\R^{n+1}$ such that
\[ 0\leq x_0 \leq x_0 ' \Rightarrow  N ( x_0 , x_1 ,\dots , x_{n}) \leq
N ( x_0 ' , x_1 ,\dots , x_{n}) \]
and
\[N(1, x_1, \dots , x_n)= 1 + \Phi (x_1 , \dots ,x_n) \]
then $\ell_\Phi = \tilde{\Lambda}(N)$.
\end{Lem}
{\bf Proof.} Let $(x_k ^1,\dots ,x_k ^n)_k \in \tilde{\Lambda} (N)$ with
$\|(x_k ^1,\dots ,x_k ^n)_k \|_{\tilde{\Lambda} (N)}\leq 1$. Let h be
the first index such that $(x_h ^1 , \dots , x_h ^n)\neq 0$. For $k
\geq h+1$ we have:
\begin{eqnarray*}
& &\underbrace{N\ast \dots \ast N}_{k}(x_1 ^1, \dots ,x_1 ^n, \dots , x_k
^1 , \dots , x_k ^n )  \\
&=& N \ast \dots \ast N(x_1 ^1 , \dots , x_{k-1} ^n)(1 + \Phi
(\frac{x_k ^1}{N \ast \dots \ast N (x_1 ^1 , \dots
,x_{k-1}^n)}, \dots ,\\
& & \dots , \frac{x_k ^n}{N \ast \dots \ast N (x_1 ^1 , \dots
,x_{k-1}^n)}))   \\
&\geq & N \ast \dots \ast N(x_1 ^1 , \dots , x_{k-1} ^n)(1 + \Phi
(x_k^1 , \dots , x_k^n))
\end{eqnarray*}
The last inequality holds because $\Phi$ is increasing on each ray
starting from $0$ and
$\|(x_k ^1,\dots ,x_k ^n)_k \|_{\tilde{\Lambda} (N)}\leq 1$. Thus
\[ \prod_{k=1}^{\infty} ( \Phi(x_k^1 , \dots , x_k^n) +1 ) < \infty
\]
and hence $\sum_{k=1}^{\infty} \Phi(x_k^1 , \dots , x_k^n) < \infty$,
 i.e. $(x_k^1 , \dots , x_k^n)_k \in \ell_\Phi$.

\smallskip
\noindent
Conversely if $(x_k^1 , \dots , x_k^n)_k \not \in \tilde{\Lambda}(N)$.
Then there exists $h$ such that\\
 $\underbrace{N\ast \dots \ast
N}_{h}(x_1 ^1 , \dots , x_h ^n)>1$. By a similar argument we see that
\[ \prod_{k=h+1}^{\infty} ( \Phi(x_k^1 , \dots , x_k^n) +1 ) = \infty
\]
which concludes the proof. \hfill $\Box$

\noindent
The next proposition \ref{P:smooth} and lemmas \ref{L:cond} and
\ref{L:norm} show how the conditions of
lemma \ref{L:suff}   can be satisfied. Let $B(x,r)$
denote the ball in $\R^n$ (with the Euclidean norm $\|\cdot\|_2$ )
centered at $x$ with radius $r$.
\begin{Prop} \label{P:smooth}
If $\phi : \R^n \rightarrow \R_+$ is a Young's function there exists\\
$\Phi: B(0,1) \rightarrow \R_+$ convex, even, $C^1$ on $B(0,1)\setminus
\{0\}$\\
 with $\Phi\sim \phi$ on $B(0,1)$.
\end{Prop}
{\bf Proof.} Let us begin with
\begin{Obs}
If $\phi : \R^n \rightarrow \R_+$ is a Young's function there exists
$\tilde{\phi}: \R^n \rightarrow \R_+$ continuous, equal to $\phi$ on
$B(0,1)$  such that
\begin{eqnarray} \label{growth}
 \lim_{\|x\|_2 \rightarrow \infty} \frac{\tilde{\phi}(x)}{\|x\|_2} =
\infty.
\end{eqnarray}
\end{Obs}
Indeed, just let
\[ \tilde{\phi}(x)= \left\{ \begin{array}{ll}
				\phi (x) & \mbox{, if } \|x\|_2 \leq 1
				\\
				\phi ( \frac{x}{ \|x\|_2 }) + (
				\|x\|_2 -1 )^2  &\mbox{, if } \|x\|_2
				> 1
			    \end{array}
		   \right.     \]

\noindent
The proof of the proposition will follow from the next two lemmas.
\begin{Lem} \label{L:smooth}
If  $\phi : \R^n \rightarrow \R_+$ is continuous and
$\phi (x)=0 \Leftrightarrow x=0$\\
then there exists $\tilde{\phi} : \R^n \rightarrow \R_+$, which is
$C^1$ on $\R^n \setminus \{0\}$ with\\
$\frac{1}{2} \phi \leq \tilde{\phi} \leq 2 \phi$.
\end{Lem}
{\bf Proof.} As $\phi$ is continuous
\[ \frac{1}{|B(x,r)|} \int_{B(x,r)} \phi \rightarrow \phi (x)
\mbox{ as } r \rightarrow 0. \]
Hence $\forall x \in \R^n \setminus \{ 0 \}$, there exists $r(x)> 0$ such
that
\begin{itemize}
\item {for $0<r<r(x)$  we have   $\frac{1}{2} \phi \leq
\frac{1}{|B(x,r)|} \int_{B(x,r)} \phi \leq 2 \phi$ and}
\item{ the map $ x \mapsto r(x)$ is continuous.}
\end{itemize}
Moreover there exists a function $\tilde{r}$, which is $C^1$ on $\R^n
\setminus \{ 0 \} $, such that $0< \tilde{r}(x) \leq r(x)$. Indeed, if
\[f(x)=\min \{r(y), y\in \overline{B}(0,\frac{1}{2^n})\setminus
B(0,\frac{1}{2^{n+1}})\} \mbox{ for } x \in
\overline{B}(0,\frac{1}{2^n})\setminus
B(0,\frac{1}{2^{n+1}}) \]
then $f_1$, the restriction of f to the positive $x_1$ axis, is a
positive step
function and we can easily choose a $C^1$ function $g$ on the positive
$x_1$ axis such that $0<g \leq f_1$. Then the radial extension of $g$
gives such an $\tilde{r}$. \\
Define
\[ \tilde{\phi}(x) = \left\{ \begin{array}{ll}
				\frac{1}{|B(x,\tilde{r}(x))|}
				\int_{B(x,\tilde{r}(x))} \phi(y)dy &
				\mbox{ if } x\neq 0 \\
				0		&  \mbox{ if } x=0
			     \end{array}
		    \right. \]
Clearly $\tilde{\phi}$ satisfies the properties required in the
conclusion of the lemma. \hfill $\Box$

\noindent
Note also that if in the previous lemma $\phi$ is even, so is
$\tilde{\phi}$. Moreover $\tilde{\phi}$ satisfies the growth condition
(\ref{growth}), if $\phi$ does.

\noindent
The next lemma follows the idea of Corollary 3.1 in \cite{GR}
(for a similar result in the infinite dimensional case see \cite{CB}).
\begin{Lem} \label{L:co}
Let $p\in \R^n$.
Let $f: \R^n \rightarrow \R$ be differentiable on $\R^n \setminus \{ p \}
$, continuous at $p$, with $ \lim_{\|x\|_2 \rightarrow \infty}
\frac{f(x)}{\|x\|_2} = \infty$. Then $\co f$ is $C^1$ on  $\R^n
\setminus \{ p \}$ .
\end{Lem}
{\bf Proof.} Let $x\in \R^n \setminus \{ p \}$. By theorem 2.1 in
 \cite{GR} there exist $q \leq n+1$, $\lambda_1 , \dots , \lambda_q$
$x_1, \dots, x_q \in \R^n$ such that:
\[\co f(x) = \sum_{i=1}^{q} \lambda_i f(x_i) \mbox{ with
}\sum_{i=1}^{q} \lambda_i x_i = x \mbox{ and }  \sum_{i=1}^{q}
\lambda_i =1 \]
As $x \neq p$ we may assume, without loss of generality, that $x_1
\neq p$. Let $U_1$ be a small neighborhood of $x_1$. For $x_1^{'} \in
U_1$ consider $x' = \lambda_1 x_1 ' + \sum_{i=2}^{n+1} \lambda_i
x_i$. Then $U=\{x'| x_1 ' \in U_1 \}$ is a neighborhood of
$x$. We can choose $U_1$ small enough such that $p \not\in U$ and $p
\not\in U_1$. Let $h: U \rightarrow U_1$ by
$h(x')=x_1 '$. Clearly h is $C^1$.
\begin{eqnarray*} \label{subgrad}
\co f(x') \leq  \lambda_1 f(x_1 ')+ \sum_{i=2}^{q} \lambda_i
f(x_i) = \lambda_1 f(h(x '))+ \sum_{i=2}^{q} \lambda_i
f(x_i)
\end{eqnarray*}
\noindent
The right hand side is a $C^1$ function of $x'$ on $U$, call it
$s(x')$. Note that $s(x)= \co f(x)$. Hence:
\begin{eqnarray}
s(y)-s(x) \geq \co f(y) -  \co f(x)
,\mbox{ for all } y \in U
\end{eqnarray}
Recall that for a convex function $\psi$ on $\R^n$ the subdifferential
of $\psi$ is a map $\partial \psi : \R^n \rightarrow \cal{P} (\R^n)$ given
by $x^{*} \in \partial \psi (x)$ if $ \psi (z) \geq \psi(x) + \langle x^* ,
z - x \rangle$ for all $z$.
\noindent
As $x$ is in the interior of the domain of $\co f$ we
have that $\partial \co f(x)$ is nonempty (see \cite{R},
Theorem 23.4). Note that if $x^* \in \partial \co f(x)$
then (\ref{subgrad}) shows that $x^* = \nabla s(x)$. Hence $\partial\co
f(x)$ is a singleton and therefore $\co
f(x)$ is differentiable at $x$ (see \cite{R}, Theorem
25.1). Hence $\co f$ is differentiable on $\R^n \setminus
\{p\}$. The conclusion of the lemma follows once we notice that if a
finite convex function on $\R^n$ is differentiable on a set then its
gradient is continuous on that set (see \cite{R}, Theorem 25.5).
\hfill $\Box$

\noindent
To prove Proposition \ref{P:smooth} let $\Phi = ( \co
\tilde{\phi})|_{B(0,1)}$ be the restriction of $\co\tilde{\phi}$ on
${B(0,1)}$, with $\tilde{\phi}$ given by lemma
\ref{L:smooth} satisfying growth condition (\ref{growth}):
smoothness of $\Phi$ is provided by lemma (\ref{L:co})
and equivalence to $\phi$ on the unit ball is obvious.
\hfill $\Box$
\begin{Lem} \label{L:cond}
Let $\phi : B(0,1) \rightarrow \R_+$ be convex, even and $C^1$ on
$B(0,1) \setminus \{0 \}$ such that $\phi (x)=0$ if and only if $x=0$.
Then there exists a Young's function $\tilde{\phi} : \R^n \rightarrow
\R_+$, which coincides with $\phi$ on a neighborhood of $0$, such that
\begin{itemize}
\item{$\forall x \in \R^n$ the map $t \mapsto \frac{1+
\tilde{\phi}(tx)}{t}$ is decreasing on $(0,\infty)$ and}
\item{$ x \mapsto \lim_{t \rightarrow \infty} \frac{1+
\tilde{\phi}(tx)}{t}$ is a norm on $\R^n$.}
\end{itemize}
\end{Lem}
{\bf Proof.} For all $\alpha > 0$ small enough $\phi^{-1} (\{ \alpha
\})$ is an $n-1$ dimensional closed submanifold of $B(0,1)$. Such an
$\alpha$ will be chosen later on. Let $| \cdot |_\alpha$ be the Minkowski
norm of the set $\phi^{-1} ([0, \alpha ])$. Then $| \cdot |_\alpha$ is
$C^1$ on $\R^n \setminus \{ 0 \}$ since $\phi$ is. An easy calculation
shows that
\begin{eqnarray} \label{nabla}
\nabla | \cdot |_\alpha
(x) = \frac{1}{\langle \nabla \phi (x), x\rangle} \nabla \phi
(x) ,\mbox{ for all } x \in \R^n \mbox{ with } |x |_\alpha =1
\end{eqnarray}
where $\langle \cdot ,\cdot \rangle$ denotes the Euclidean inner
product on $\R^n$.
Indeed, for $x \in \R^n \setminus \{ 0 \}$ , $ |x|_\alpha =
\frac{1}{\lambda (x)}$ where $ \phi( \lambda(x) x) = \alpha$.
Thus
\[ \nabla  | \cdot |_\alpha (x) = - \frac{1}{\lambda ^2 (x)} \nabla \lambda
(x)= - | x |_\alpha^2 \nabla \lambda(x) \]
By differentiating $\phi(\lambda(x)x) = \alpha$ with respect to $x_j$
we obtain
\begin{eqnarray*}
0&=&\frac{\partial}{\partial x_j} \phi (\lambda(x) x) = \sum_{i=1}^n
\frac{ \partial \phi}{\partial x_i}(\lambda(x) x) \left( \frac{\partial
\lambda}{\partial x_j}(x) x_i + \lambda(x) \frac{\partial
x_i}{\partial x_j}\right) \\
&=&\frac{\partial \phi}{\partial x_j}(\lambda(x)x) \lambda(x) +
\sum_{i=1}^n \frac{\partial \phi}{\partial x_i}(\lambda(x)x)
\frac{\partial \lambda}{\partial x_j}(x) x_i \\
&=&\frac{\partial \phi}{\partial x_j}\left( \frac{x}{|x|_\alpha}\right)
\frac{1}{|x|_\alpha}
+ \langle \nabla \phi \left( \frac{x}{|x|_\alpha} \right), x \rangle
\frac{\partial
\lambda}{\partial x_j}(x)
\end{eqnarray*}
Thus
\[\frac{\partial \lambda}{\partial x_j}(x) = - \mbox{ }\frac{1}{|x|_\alpha
\langle \nabla
\phi \left( \frac{x}{|x|_\alpha} \right), x \rangle } \mbox{ } \frac{\partial
\phi}{\partial x_j} \left( \frac{x}{|x|_\alpha} \right) \]
and therefore
\begin{eqnarray*}
\nabla | \cdot |_\alpha (x) &=&- |x|_\alpha^2 \left( \frac{-1}{|x|_\alpha
\langle \nabla
\phi  \left( \frac{x}{|x|_\alpha}\right), x \rangle} \right) \nabla
\phi \left(
\frac{x}{|x|_\alpha} \right)\\
&=&  \frac{|x|_\alpha}{\langle \nabla \phi \left(
\frac{x}{|x|_\alpha} \right), x\rangle} \nabla \phi \left(
\frac{x}{|x|_\alpha}
\right)
\end{eqnarray*}
which gives (\ref{nabla}) for $|x|_\alpha=1$.
Let $M=\sup_{|x|_\alpha=1} \langle \nabla \phi (x) , x \rangle$. Then
for every $x$
with $|x|_\alpha =1$ and for every $ u\in \R^n$ such that $\langle \nabla |
\cdot |_\alpha (x), u \rangle > 0$ we have
\[ \langle \nabla \phi (x) , u \rangle \leq M \langle \nabla | \cdot
|_\alpha (x), u \rangle \]
Thus, for all such $x$ and $u$ we have
\begin{eqnarray} \label{dir}
D_u \phi (x) \leq D_u (M | \cdot |_\alpha ) (x)
\end{eqnarray}
where $D_u$ is the directional derivative in the direction of $u$.

\noindent
Define $\tilde{\phi}$ on $\R^n$ by:
\[ \tilde{\phi}(x) = \left\{  \begin{array}{lll}
				\phi (x)         & \mbox{, if $|x|_\alpha
				\leq 1$}\\
				\alpha + M ( |x|_\alpha -1) &
				\mbox{, otherwise}
			      \end{array}
		    \right. \]

\noindent
Condition (\ref{dir}) provides the convexity of $\tilde{\phi}$.
Clearly $\tilde{\phi}$ is a Young's function and coincides with $\phi$
in a neighborhood of $0$.

\noindent
Fix $x \in \R^n \setminus \{ 0 \}$. We want to show that the mapping
\[ t\mapsto \left\{ \begin{array}{ll}
			\frac{1+ \phi(tx)}{t} & \mbox{, if } 0< t \leq
			\frac{1}{|x|_\alpha}\\
			\frac{1+ \alpha + M ( |tx|_\alpha -1)}{t} & \mbox{,
			if } t \geq \frac{1}{|x|_\alpha}
		   \end{array}
	   \right. \]

\noindent
is decreasing.

\noindent
As the map is continuous, it suffices to show it is decreasing on $
(0, \frac{1}{|x|_\alpha})$  and on $ (\frac{1}{|x|_\alpha} , \infty)$. Since
$\phi$ is convex and $\phi(0)=0$ there exists $\tau=\tau(x)$ such that
$ t\mapsto \frac{1+\phi(tx)}{t}$  is decreasing on  $(0, \tau(x)
)$ with $0< \tau(x) \leq \frac{1}{\|x\|_2}$ and $\tau(\mu
x)=\frac{\tau(x)}{\mu}$.
\noindent
By compactness of $\overline{B}(0, 1)$ and the continuity of $\tau$
and $\phi$ we have $\inf_{\|y\|_2 =\frac{1}{2}} \phi ( \tau(y)y) >
0$. If
\begin{eqnarray} \label{alfa1}
\alpha \leq \inf_{\|y\|_2 =\frac{1}{2}} \phi ( \tau(y)y)
\end{eqnarray}
then
\[ \phi \left( \tau \left( \frac{x}{2 \|x\|_2} \right) \frac{x}{2
\|x\|_2} \right) > \alpha \]
hence $\phi( \tau (x)x) > \alpha$ and thus $| \tau(x) x |_\alpha >1$.
Therefore $\frac{1}{|x|_\alpha}< \tau(x)$ and $ t\mapsto
\frac{1+\phi(tx)}{t}$  is decreasing on $(0,\frac{1}{|x|_\alpha})$.

\noindent
Note that $\frac{1+ \alpha + M ( t |x|_\alpha -1)}{t} = \frac{1+ \alpha -
M}{t} + M |x|_\alpha$ is decreasing as a function of $t$ exactly when $1+
\alpha - M > 0$. Thus it suffices to have $M\leq 1$. But
\begin{eqnarray*}
M&=& \sup_{|x|_\alpha=1} \langle \nabla \phi (x) , x \rangle \mbox{ }
= \sup_{|x|_\alpha=1} \langle \nabla \phi (x),\frac{x}{\|x\|_2}\rangle
\|x\|_2 \\
&\leq& \left(\sup_{|x|_\alpha=1}
D_{\frac{x}{\|x\|_2}}\phi(x)\right)(\sup_{|x|_\alpha=1}\|x\|_2)
\rightarrow 0 \mbox{ as } \alpha \rightarrow 0
\end{eqnarray*}
since $D_{\frac{x}{\|x\|_2}}\phi(x)$ is an increasing function of $\alpha$
and $\sup_{|x|_\alpha=1} \|x\|_2 \rightarrow 0$ as $\alpha \rightarrow 0$.
In particular $\alpha$ can be chosen such that (\ref{alfa1}) holds and
$M \leq 1$. Then $t \mapsto \frac{1+ \tilde{\phi}(tx)}{t}$ is
decreasing $\forall x\in \R^n$.

\noindent
Finally, note that
\[x \mapsto \lim_{t \rightarrow \infty} \frac{1+ \tilde{\phi}(tx)}{t}
= M |x|_\alpha \]
is a norm on $\R^n$. \hfill $\Box$

\begin{Lem} \label{L:norm}
Let $\phi: \R^n \rightarrow \R_+$ be a Young's function such that
\begin{itemize}
\item{$\forall x \in \R^n$ the map $t \mapsto \frac{1+
\phi(tx)}{t}$ is decreasing on $(0,\infty)$}
\item{$ x \mapsto \lim_{t \rightarrow \infty} \frac{1+
\phi(tx)}{t}$ is a norm on $\R^n$}
\end{itemize}
Then
\[ N(x_0 , x_1 , \dots ,x_n ) = \left\{ \begin{array}{ll}
					  |x_0 | (1 + \phi
					  \left(\frac{ x_1 , \dots
					  ,x_n }{|x_0 |}\right)) \mbox{ if }
					  x_0 \neq 0 \\
					  \lim_{t \rightarrow 0}|t | (1 + \phi
					  \left(\frac{ x_1 , \dots
					  ,x_n }{|t|}\right))
					  \mbox{ if } x_0 =0
					\end{array}
					\right. \]

\noindent
is a norm on $\R^{n+1}$ such that
\[ 0\leq x_0 \leq x_0 ' \Rightarrow  N ( x_0 , x_1 ,\dots , x_{n}) \leq
N ( x_0 ' , x_1 ,\dots , x_{n}) \]
and
\[N(1, x_1, \dots , x_n)= 1 + \Phi (x_1 , \dots ,x_n). \]
\end{Lem}
{\bf Proof.}  We only need to check the triangle inequality for
$N$. Let $x= (x_0 , x_1 , \dots , x_n)$ and $ y = (y_0 , y_1 , \dots ,
y_n) \in \R^{n+1}$.

\noindent
If $x_0 > 0$ and $y_0 > 0$
\begin{eqnarray*}
& &N(x+y) = (x_0 + y_0 )\left( 1 + \phi \left( \frac{x_1 + y_1}{x_0 + y_0},
\dots ,\frac{x_n + y_n}{x_0 + y_0}\right)\right) \\
&=& x_0 + y_0 + (x_0 + y_0 )\phi \left(\frac{x_0}{x_0 + y_0} \frac{x_1
, \dots , x_n}{x_0} + \frac{y_0}{x_0 + y_0} \frac{y_1 , \dots ,
y_n}{y_0} \right) \\
&\leq & x_0 + y_0 + x_0 \phi \left( \frac{x_1, \dots , x_n}{x_0}
\right) + y_0 \phi \left(\frac{y_1 , \dots ,y_n}{y_0} \right) \\
&=& N(x) + N(y)
\end{eqnarray*}
by the convexity of $\phi$.

\noindent
If  $x_0 > 0$ and $y_0 = 0$
\begin{eqnarray*}
N(x+y)  \leq  N(x_0 - \epsilon ,x_1, \dots, x_n) + N(\epsilon , y_1,
\dots, y_n)  \mbox{, for all } \epsilon \in (0, x_0)
\end{eqnarray*}
by the previous case. Letting $\epsilon \rightarrow 0$, by the
continuity of $N$ we obtain the desired inequality.

\noindent
Finally, if  $x_0 > 0$ and $y_0 < 0$ we may assume $0 \leq x_0 + y_0 <
x_0$. Then, by the properties of $\phi$ and the previous case we have
\begin{eqnarray*}
N(x+y) & \leq & N(x_0 , x_1 + y_1 , \dots ,x_n +y_n) \\
& \leq & N(x) + N(0, y_1 , \dots , y_n) \\
& \leq & N(x) + N(|y_0 | , y_1 , \dots , y_n) = N(x) + N(y)
\end{eqnarray*}
which concludes the proof. \hfill $\Box$

\smallskip
\noindent
The proof of theorem \ref{T:FOproM} is now complete.

\smallskip
\noindent
Theorems \ref{T:twistOFO} and \ref{T:FOproM} give immediately the
following
\noindent
\begin{Cor}
Let $\ell_\phi$ be an Orlicz space with non-trivial type. Then the twisted sum
$\ell_\phi \bigoplus_F \ell_\phi$ can be equivalently renormed to have
property (M).
\end{Cor}
In particular, the spaces $Z_p , 1<p<\infty,$ have property (M) after
renorming. We end
this section with an application of the previous corollary and Theorem
2.4 in \cite{K4}. Recall that if $X$ is a Banach space and $E$ is a
subspace of $X$ then $E$ is called an {\it M-ideal in $X$} (see \cite{AE}) if
$X^*$ can be decomposed as an $\ell_1$-sum $X^* =E^{\perp}\bigoplus_1
V$ for some closed subspace $V$ of $X^*$. For a Banach space $X$ let
${\cal{L}}(X)$ denote the algebra of all bounded operators on $X$ and
${\cal{K}}(X)$ the ideal of compact operators.
\begin{Cor}
Let $\ell_\phi$ be an Orlicz space with non-trivial type. Then\\
${\cal{K}}(\ell_\phi \bigoplus_F \ell_\phi)$ is an M-ideal in
${\cal{L}}( \ell_\phi \bigoplus_F \ell_\phi).$
\end{Cor}

\end{document}